# MEANS OF A DIRICHLET PROCESS AND MULTIPLE HYPERGEOMETRIC FUNCTIONS[1]

By Antonio Lijoi and Eugenio Regazzini[2]

*Università degli Studi di Pavia*

The Lauricella theory of multiple hypergeometric functions is used to shed some light on certain distributional properties of the mean of a Dirichlet process. This approach leads to several results, which are illustrated here. Among these are a new and more direct procedure for determining the exact form of the distribution of the mean, a correspondence between the distribution of the mean and the parameter of a Dirichlet process, a characterization of the family of Cauchy distributions as the set of the fixed points of this correspondence, and an extension of the Markov–Krein identity. Moreover, an expression of the characteristic function of the mean of a Dirichlet process is obtained by resorting to an integral representation of a confluent form of the fourth Lauricella function. This expression is then employed to prove that the distribution of the mean of a Dirichlet process is symmetric if and only if the parameter of the process is symmetric, and to provide a new expression of the moment generating function of the variance of a Dirichlet process.

**1. Introduction.** The connections between properties of functionals of a Dirichlet process and the Lauricella multiple hypergeometric functions have been investigated in independent papers by Kerov and Tsilevich (1998) and Regazzini (1998); they also represent the point of departure of the present paper. The approach undertaken here is quite different from that of recent contributions to the study of exact distributions of functionals of a Dirichlet process. See, for example, Regazzini, Guglielmi and Di Nunno (2002). Some

Received April 2002; revised May 2003.

[1]Supported in part by MURST, Programma di Ricerca "*Metodi bayesiani non parametrici e loro applicazioni,*" 2002–2004, and by Università di Pavia, Progetto di Ricerca di Ateneo "*Analisi statistica di sistemi complessi.*"

[2]Also affiliated with IMATI-CNR, Milano, Italy.

specific properties of multiple hypergeometric functions are extended in such a way as to become significant properties of the laws of functionals of random measures. This is an unusual application of the theory of special functions. On the other hand, since these extensions can be thought of as infinite-dimensional versions of some fundamental types of special functions, the results presented in the following sections might be also of some interest from an analytic point of view.

The main reason for analyzing the interplay between multiple hypergeometric functions and laws of Dirichlet functionals is twofold: it allows simplifications in the proofs of some well-known propositions and, more importantly, it leads to new results concerning the distribution of the abovementioned functionals. Ties between the theory of multiple hypergeometric functions and some applied problems in statistics had been investigated in previous contributions, in addition to the papers already cited above. Dickey (1968) extends an identity due to Picard and applies it to Bayesian inference about multinomial cell probabilities, with a prior expressed by the Savage generalization of the Dirichlet distribution. Hill (1977) exploits suitable representations of Appell's hypergeometric functions to exactly evaluate the posterior moments of parameters of interest for the inference about variance components. In a problem of Bayesian statistical inference for missing data, Dickey, Jiang and Kadane (1987) obtain representations of posterior moments and predictive probabilities in terms of ratios of Carlson's hypergeometric functions. Jiang, Kadane and Dickey (1992) explore computational methods for hypergeometric functions arising in Bayesian analysis. Moreover, an introduction to Carlson's functions for statisticians can be found in Dickey (1983).

As far as the present article is concerned, it is organized as follows. Some integral representations of the Lauricella hypergeometric functions, together with their probabilistic interpretations, are illustrated in Section 2. The Feigin and Tweedie existence condition and the Markov–Krein identity are then jointly restated in Theorem 1 in Section 3. The proof is based on some classical results concerning Lauricella's hypergeometric functions. Despite its pure analytic nature, it is as simple as the proof of Theorem 9 in Tsilevich, Vershik and Yor (2000). See also Theorem 2 in Tsilevich, Vershik and Yor (2001). Moreover, we prove the existence of a one-to-one correspondence between the set of nonnull measures, on $\mathbb{R}$, with given finite total mass $a > 0$ and the set of all probability distributions (p.d.) of linear functionals of a Dirichlet process with parameter having total mass $a$. Among recent papers focusing on such a bijection—sometimes referred to as the *Markov–Krein correspondence*—are those of Diaconis and Kemperman (1996), Kerov and Tsilevich (1998) and Tsilevich, Vershik and Yor (2000). In Section 4, elementary properties of gamma processes lead to an expression for the p.d. of the mean of a gamma process in terms of the p.d. of the mean of a Dirichlet process. Moreover, it is also seen that the p.d. of the mean of a gamma process is a generalized gamma convolution and is, hence, infinitely divisible. This fact suggests a simple proof of the absolute continuity, with respect to the Lebesgue measure on $\mathbb{R}$, of the mean of a Dirichlet process (see Proposition 2 and Remark 1). In Section 5, an extension of the Markov–Krein correspondence is deduced from well-known integral representations of the fourth Lauricella hypergeometric function. Two applications of this extension are considered in Section 6: the first determines the exact form of the p.d. of the mean of a Dirichlet process, and the second is a characterization of the Cauchy distribution. Section 7 is devoted to a representation of the characteristic function of the mean of a Dirichlet process, via a multidimensional extension of Kummer's confluent hypergeometric function. The latter can be considered as a confluent form of the fourth Lauricella function and admits a representation as a contour integral proved by Erdélyi (1937). The distribution of a vector of means of a single Dirichlet process is examined in Section 8. The identities given in Section 5 are trivially extended to this setup. This is also helpful in proving the absolute continuity, with respect to the Lebesgue measure on $\mathbb{R}^d$, of a vector of means of a single Dirichlet process, a finding which answers a question raised in Regazzini, Guglielmi and Di Nunno (2002). Finally, a representation of the moment generating function of the variance of the Dirichlet process is provided and it would seem a meaningful improvement on previous contributions to the subject. See, for example, Cifarelli and Melilli (2000).

**2. Probabilistic interpretation of the Lauricella fourth function $F_D$.** The topic of multiple hypergeometric functions was first approached, in a systematic way, by Giuseppe Lauricella at the end of the 19th century. See, for example, Exton (1976). He defined four functions which are named after him and have both multiple series and integral representations. In particular, the *fourth* of these functions, denoted by $F_D$, admits integral representations of importance in our treatment. Given any $\boldsymbol{\xi} = (\xi_1, \ldots, \xi_n)$ in $\mathbb{R}^n$, set $|\boldsymbol{\xi}|$ for $\sum_{k=1}^n \xi_k$ and $\langle \cdot, \cdot \rangle$ for inner product. Moreover, let $T_n := \{\mathbf{u} = (u_1, \ldots, u_n) \in \mathbb{R}^n : u_i > 0 \text{ for every } i \text{ and } |\mathbf{u}| < 1\}$. With this notation, an integral representation of Euler type [see (26) in Lauricella (1893)] is

$$
\frac{\Gamma(b_1) \cdots \Gamma(b_n) \Gamma(a - |\mathbf{b}|)}{\Gamma(a)} F_D(c, b_1, \ldots, b_n; a; x_1, \ldots, x_n)
$$

(2.1)

$$
= \int_{T_n} u_1^{b_1 - 1} \cdots u_n^{b_n - 1} (1 - |\mathbf{u}|)^{a - |\mathbf{b}| - 1} (1 - \langle \mathbf{u}, \mathbf{x} \rangle)^{-c} \, du_1 \cdots du_n,
$$

which is valid for every $\mathbf{x} = (x_1, \ldots, x_n)$ in $[0, 1)^n$ and $\mathbf{b} = (b_1, \ldots, b_n)$ with strictly positive real part, that is, $\operatorname{Re}(b_k) > 0$ for every $k$, and such that $\mathrm{Re}(a - |\mathbf{b}|) > 0$. A further representation provided by formula (25) in Lauricella ([1893](#)) is

$$
\begin{aligned}
(2.2) \quad & F_D(c, b_1, \ldots, b_n; a; x_1, \ldots, x_n) \\
& = \int_{[0,1]} (1 - u x_1)^{-b_1} \cdots (1 - u x_n)^{-b_n} B(\, du; c, a - c),
\end{aligned}
$$

which holds true whenever $\mathrm{Re}(c) > 0$ and $\mathrm{Re}(a - c) \geq 0$, with

$$
B(A; c, a - c) = \begin{cases}
\delta_1(A), & \text{if } \mathrm{Re}(a - c) = 0, \\[2mm]
\dfrac{\Gamma(a)}{\Gamma(c)\Gamma(a - c)} & \\[1mm]
\quad \times \int_{A \cap (0,1)} u^{c-1} (1 - u)^{a-c-1} \, du, & \text{if } \mathrm{Re}(a - c) > 0,
\end{cases}
$$

for every $A$ in the Borel class on $\mathbb{R}$, $\mathscr{B}(\mathbb{R})$, and $\delta_x$ standing for the unit point mass concentrated at $x$. Hence, letting $a = c$, one gets

$$
(2.3) \quad F_D(a, b_1, \ldots, b_n; a; x_1, \ldots, x_n) = (1 - x_1)^{-b_1} \cdots (1 - x_n)^{-b_n}.
$$

Representation (2.1) has an obvious probabilistic interpretation, which was first stressed and exploited by Carlson ([1977](#)), giving rise to what he called "method of Dirichlet averages," for unifying a part of the theory of special functions. In point of fact, let $\tilde{\varphi}$ be a random probability measure supported by $S = \{0, x_1, \ldots, x_n\}$, with $(x_1, \ldots, x_n)$ in $(0,1)^n$, and assume that the random vector $(\tilde{\varphi}\{x_1\}, \ldots, \tilde{\varphi}\{x_n\})$ has the Dirichlet distribution, $\mathscr{D}$, with parameter $(b_1, \ldots, b_n, a - |\mathbf{b}|)$, and $a - |\mathbf{b}| > 0$, $b_k > 0$ for every $k$, that is, the distribution characterized by

$$
\begin{aligned}
\mathscr{D}(\, d\mathbf{u}) = {} & \frac{\Gamma(a)}{\Gamma(b_1) \cdots \Gamma(b_n)\Gamma(a - |\mathbf{b}|)} \\
& \times u_1^{b_1 - 1} \cdots u_n^{b_n - 1} (1 - |\mathbf{u}|)^{a - |\mathbf{b}| - 1} \mathbb{I}_{T_n}(\mathbf{u}) \, du_1 \cdots du_n,
\end{aligned}
$$

where $\mathbb{I}_B$ denotes the indicator function of set $B$. Whence, combination of (2.1) with (2.3) gives

$$
(2.4) \quad \mathscr{D}\left[ \left( 1 - \int_S x \tilde{\varphi}(dx) \right)^{-a} \right] = \exp\left\{ - \int_S \log(1 - x) \alpha(dx) \right\},
$$

provided that $\alpha$ is the measure on the power set of $S$ determined by $\alpha\{x_k\} = b_k$ for $k = 1, \ldots, n$ and $\alpha\{0\} = a - |\mathbf{b}|$. Here, and in the sequel, given any probability measure $\varphi$ and real-valued function $h$ such that $\int |h| \, d\varphi < +\infty$, $\varphi(h)$ denotes $\int h \, d\varphi$.

Equality (2.4) represents the most elementary version of an identity established by Cifarelli and Regazzini ([1979a](#), [b](#), [1990](#)) for *functional Dirichlet processes* with parameter $\alpha$. Recall that, given a finite measure $\alpha$ on $(\mathbb{R}, \mathscr{B}(\mathbb{R}))$ with $a := \alpha(\mathbb{R}) > 0$, a random probability measure $\tilde{\varphi}$ is said to be a (functional) Dirichlet process with parameter $\alpha$ if, for every finite and measurable partition $C_1, \ldots, C_n, C_{n+1}$ of $\mathbb{R}$, the random vector $(\tilde{\varphi}(C_1), \ldots, \tilde{\varphi}(C_n))$ has the Dirichlet distribution with parameter $(\alpha(C_1), \ldots, \alpha(C_n), \alpha(C_{n+1}))$. As for the definition and the main properties of a functional Dirichlet process, the seminal contribution in Ferguson ([1973](#)) still represents a sound reference.

Throughout the following sections, $\zeta(w; \alpha, f)$ will indicate the integral $\int_{\mathbb{R}} \log(1 + wf(x)) \alpha(dx)$, while $\log z$ will denote the principal determination of the logarithm of the complex number $z$, that is, $\log z = \log|z| + i \arg(z)$, where $\arg(z)$ is chosen in $(-\pi, \pi]$.

## 3. Lauricella theory and Markov–Krein identity.

Diaconis and Kemperman observed that the identity of Cifarelli and Regazzini cited above is closely related to a well-known version of the Markov moment problem. See Akhiezer and Krein ([1962](#)), Diaconis and Kemperman ([1996](#)) and Kerov ([1998](#)). As a matter of fact, some authors refer to this identity as the *Markov–Krein identity*. Compare, for example, Kerov and Tsilevich ([1998](#)) and the more recent paper by Tsilevich, Vershik and Yor ([2000](#)). According to Kerov and Tsilevich ([1998](#)) and Regazzini ([1998](#)), the Markov–Krein identity extends ([2.4](#)) from the (finite-dimensional) Dirichlet distribution to the functional Dirichlet distribution. Notice that in the real domain the extension holds true provided that the support of the parameter of the latter distribution is bounded above. In the complex domain, a variant of the same identity, which holds for any parameter, is established in Kerov and Tsilevich ([1998](#)) and in Regazzini, Guglielmi and Di Nunno ([2002](#)). In the following section, the Lauricella theory is exploited to state jointly the Markov–Krein identity and the Feigin and Tweedie condition for finiteness of the mean of a Dirichlet process.

### 3.1. *The Markov–Krein identity.*

Given any finite measure $\alpha$ on $(\mathbb{R}, \mathscr{B}(\mathbb{R}))$ such that $\alpha(\mathbb{R}) = a > 0$, let $\mathscr{D}_\alpha$ denote the functional Dirichlet p.d. with parameter $\alpha$. A suitable sequence of simple functions $\xi_n$, $n \geq 1$, exists that satisfies $\xi_n(x) \uparrow I(x)$ if $x \geq 0$, $\xi_n(x) \downarrow I(x)$ when $x < 0$ and $|\xi_n| \uparrow |I|$, where $I$ is the identity map on $\mathbb{R}$. Therefore the sequence of finite measures $\alpha_n$ defined by

$$\alpha_n(B) = \alpha\{\xi_n \in B\}, \qquad B \in \mathscr{B}(\mathbb{R})$$

converges weakly to $\alpha$, namely $\alpha_n \Rightarrow \alpha$. Now define $\mathbb{P}$ to be the space of all probability measures on $(\mathbb{R}, \mathscr{B}(\mathbb{R}))$ endowed with the topology of weak convergence, and let $\mathscr{P}$ be the Borel $\sigma$-algebra on $\mathbb{P}$. In this way, it follows that the identity map $\tilde{\varphi}$ on $(\mathbb{P}, \mathscr{P}, \mathscr{D}_\alpha)$ is a random probability measure with p.d. $\mathscr{D}_\alpha$. Moreover, $\tilde{\varphi}_n := \tilde{\varphi} \circ \xi_n^{-1}$ has p.d. $\mathscr{D}_{\alpha_n}$. A straightforward application of the Lauricella formula (2.4), combined with a standard analytic continuation argument, yields

$$(3.1) \qquad \int_{\mathbb{P}} (1 + it\varphi(f \circ \xi_n))^{-a} \mathscr{D}_{\alpha}(d\varphi) = \exp\{-\zeta(it; \alpha_n, f)\}.$$

**Theorem 1.** *Setting* $L := \{\varphi \in \mathbb{P} : \varphi(|I|) := \lim \varphi(|\xi_n|) \text{ is finite}\}$, *one has*

$$\mathscr{D}_{\alpha}(L)\{1 - \mathscr{D}_{\alpha}(L)\} = 0,$$

*and the following two conditions are equivalent:*

(i) $\mathscr{D}_{\alpha}(L) = 1$.
(ii) $\int \log(1 + |x|)\alpha(dx) < +\infty$.

*Moreover, if* $\mathscr{D}_{\alpha}(L) = 1$, *then*

$$(\text{iii}) \qquad \int_{\mathbb{P}} \frac{1}{(1 + it\varphi(I))^a} \mathscr{D}_{\alpha}(d\varphi) = \exp\{-\zeta(it; \alpha, I)\}, \qquad t \in \mathbb{R}.$$

The equivalence of (i) and (ii) was first proved by Feigin and Tweedie (1989), but it was already contained—at least in part—in Cifarelli and Regazzini (1979a, b). A complete development of the argument used in the latter paper can be found in Cifarelli and Regazzini (1996). Here it is shown that both equivalence of (i) and (ii) and identity (iii) easily follow from (3.1).

Proof of Theorem 1.    In view of the definition of $(\xi_n)_{n \geq 1}$, the operations of integration and of taking limit can be interchanged to obtain

$$\zeta(it; \alpha_n, |I|)$$

$$= \tfrac{1}{2} \int_{\mathbb{R}} \log(1 + t^2 \xi_n^2(x))\alpha(dx) + i \int_{\mathbb{R}} \arg(t|\xi_n(x)|)\alpha(dx)$$

$$\to \tfrac{1}{2} \int_{\mathbb{R}} \log(1 + t^2 x^2)\alpha(dx) + i \int_{\mathbb{R}} \arg(t|x|)\alpha(dx), \qquad \text{as } n \to +\infty,$$

where an infinite real part in the above limiting expression is allowed. Moreover,

$$\int_{\mathbb{P}} (1 + it\varphi(|\xi_n|))^{-a} \mathscr{D}_{\alpha}(d\varphi)$$

$$= \int_{\mathbb{P}} \exp\left\{-\frac{a}{2} \log(1 + t^2 \varphi^2(|\xi_n|)) - ia \arg(t\varphi(|\xi_n|))\right\} \mathscr{D}_{\alpha}(d\varphi)$$

$$\to \int_{L} \frac{1}{(1 + it\varphi(|I|))^a} \mathscr{D}_{\alpha}(d\varphi), \qquad \text{as } n \to +\infty,$$

since $\varphi(|I|) = +\infty$ if $\varphi \in L^c$. Hence, by (3.1) with $f(\cdot) = |\cdot|$, one gets

$$(3.2) \qquad \int_{L} (1 + it\varphi(|I|))^{-a} \mathscr{D}_{\alpha}(d\varphi) = \exp\{-\zeta(it; \alpha, |I|)\}$$

for every $t \in \mathbb{R}$, with the proviso that the right-hand side is 0 whenever the real part of $\int \log(1 + it|x|)\alpha(dx)$ is infinite. Thus, if $\int_{\mathbb{R}} \log(1 + |x|)\alpha(dx) = +\infty$, the right-hand side of (3.2) is 0 for every $t \neq 0$ and, therefore, $\mathscr{D}_\alpha(L)$ must be 0. On the other hand, if $\int_{\mathbb{R}} \log(1 + |x|)\alpha(dx) < +\infty$, then $t \mapsto \exp\{-\int_{\mathbb{R}} \log(1 + it|x|)\alpha(dx)\}$ is continuous at $t = 0$, that is,

$$\lim_{t \to 0} \exp\left\{-\int_{\mathbb{R}} \log(1 + it|x|)\alpha(dx)\right\} = 1,$$

and, by taking the limit (as $t \to 0$) on both sides of (3.2), one obtains

$$\mathscr{D}_\alpha(L) = 1.$$

Finally, if (ii) holds, interchanging the operations of integration and of taking limit in (3.1), with $f = I$, yields (iii). $\quad\square$

Identity (iii) in Theorem 1 can be invoked to prove that the p.d. of $\tilde{\varphi}$ characterizes the parameter of $\mathscr{D}_\alpha$ within the class of all measures $\alpha$ on $(\mathbb{R}, \mathscr{B}(\mathbb{R}))$ for which $\alpha(\mathbb{R})$ has some fixed value $a > 0$. This is the essence of the Markov–Krein correspondence and it is dealt with in Section 3.2. In addition to the papers listed in Section 1, Andrea Ongaro, in a personal communication, has shown us a proof of such a correspondence based on a completely different approach.

3.2. *Uniqueness theorem.* Define $\mathbb{F}$ to be the set of all finite measures on $(\mathbb{R}, \mathscr{B}(\mathbb{R}))$, and put

$$\mathbb{F}_a = \left\{\alpha \in \mathbb{F} : \alpha(\mathbb{R}) = a \text{ and } \int_{\mathbb{R}} \log(1 + |x|)\alpha(dx) < +\infty\right\},$$

$$\mathbb{M}_a = \{\mathscr{D}_\alpha \circ \tilde{\varphi}(I)^{-1} : \alpha \in \mathbb{F}_a\}$$

for every $a > 0$. In words, $\mathbb{M}_a$ is the set of all p.d.'s of $\tilde{\varphi}(I)$ when $\tilde{\varphi}$ is a Dirichlet process with parameter varying in $\mathbb{F}_a$. Let $A$ stand for the distribution function associated to $\alpha$, and $\bar{A}$ for $(a - A)\mathbb{I}_{[0, +\infty)} - A\mathbb{I}_{(-\infty, 0)}$. Clearly, to any $\alpha$ in $\mathbb{F}_a$ there corresponds a unique $\mu_\alpha := \mathscr{D}_\alpha \circ \tilde{\varphi}(I)^{-1}$. Moreover, one has

THEOREM 2. *Any $\alpha$ in $\mathbb{F}_a$ is determined by its $\mu_\alpha$.*

PROOF. Fix $\mu_{\alpha_1}$ in $\mathbb{M}_a$ and suppose $\alpha_2 \in \mathbb{F}_a$ is such that $\mu_{\alpha_2} = \mu_{\alpha_1}$. Then

$$\exp\{-\zeta(it; \alpha_1, I)\} = \int_{\mathbb{R}} \frac{1}{(1 + itx)^a} \mu_{\alpha_1}(dx)$$

$$(3.3) \qquad\qquad = \int_{\mathbb{R}} \frac{1}{(1 + itx)^a} \mu_{\alpha_2}(dx) \qquad \text{[by hypothesis]}$$

$$= \exp\{-\zeta(it; \alpha_2, I)\} \qquad \text{[from (iii) in Theorem 1].}$$

Use integration by parts to obtain

$$\zeta(it; \alpha_j, I) = \int_{\mathbb{R}} \frac{it}{1 + itx} \bar{A}_j(x) \, dx, \qquad t \in \mathbb{R}, j = 1, 2,$$

which, combined with (3.3) and analytic continuation, yields

$$(3.4) \qquad \int_{\mathbb{R}} \frac{\bar{A}_1(x)}{z + x} \, dx = \int_{\mathbb{R}} \frac{\bar{A}_2(x)}{z + x} \, dx, \qquad z \in \mathbb{C} \text{ with } \operatorname{Im}(z) \neq 0.$$

Taking, $l(z, \alpha_j)$ for $\int_{\mathbb{R}} [\bar{A}_j(x)/(z+x)] \, dx$, $j = 1, 2$, and resorting to the Stieltjes–Perron inversion formula [see Theorem 12.10d in Henrici (1991)], one has

$$\begin{aligned}
\bar{A}_1(\xi) &= \frac{1}{2\pi i} \lim_{\varepsilon \downarrow 0} \{ l(-\xi - i\varepsilon, \alpha_1) - l(-\xi + i\varepsilon, \alpha_1) \} \\
&= \frac{1}{2\pi i} \lim_{\varepsilon \downarrow 0} \{ l(-\xi - i\varepsilon, \alpha_2) - l(-\xi + i\varepsilon, \alpha_2) \} \\
&= \bar{A}_2(\xi),
\end{aligned}$$

provided that $\xi$ is a continuity point for both $\bar{A}_1$ and $\bar{A}_2$. This suffices to conclude that $\alpha_1 = \alpha_2$, since $\alpha_i(\mathbb{R}) = a$, $i = 1, 2$. $\quad \square$

Theorem 2 states that there is a bijection $\beta_a$ of $\mathbb{F}_a$ to $\mathbb{M}_a$. Clearly, $\mathbb{M}_a \subset \mathbb{P}$ for every $a > 0$ and, more precisely, $\mathbb{M}_a \subsetneq \mathbb{P}$ for each $a > 0$. Namely, for every $a > 0$, there is some probability $\mu$ on $\mathscr{B}(\mathbb{R})$ which is not the p.d. of $\tilde{\varphi}(I)$ for any $\mathscr{D}_\alpha$ with $\alpha \in \mathbb{F}_a$. For example, take

$$\mu(dx) = \frac{1}{2\eta} \{ \mathbb{I}_{(-1-\eta, -1)}(x) + \mathbb{I}_{(1, 1+\eta)}(x) \} \, dx$$

with $\eta > 0$. It is easily seen that such a $\mu$ cannot be the p.d. of the mean of a Dirichlet process with parameter $\alpha$, for any choice of $\alpha$ in $\mathbb{F}_a$ and for every $a > 0$. One can argue by noting that the support of the p.d. of the mean of a Dirichlet process with parameter $\alpha$ must coincide with the closure of the convex hull of the support of $\alpha$. This can be proved by combining a result concerning the topological support of the Dirichlet prior given in Majumdar (1992) and the equivalence of (i) and (ii) in Theorem 1. We finally remark that general properties of $\mathbb{M}_a$ are discussed, for example, in Section 2.3 of Kerov (1998).

**4. Markov–Krein identity and means of gamma processes.** From elementary properties of the gamma process it follows that the right-hand side of (iii) in Theorem 1 is the conjugate of the characteristic function of the "mean" of a gamma process with parameter $\alpha$. See (4.1). It is easy to verify that such a characteristic function coincides with the Fourier–Stieltjes transform of a *generalized gamma convolution*, according to the terminology introduced by Olof Thorin. See, for example, the original contributions in Thorin ([1977a](#), b, [1978a](#), b) and the systematic treatment in Bondesson ([1992](#)). From one of the Thorin results, it is possible to state the infinite divisibility of the distribution of the mean of a gamma process. This fact motivates the search for the Lévy–Khintchine representation of the above characteristic function. And this representation, in turn, allows us to deduce immediately the absolute continuity, with respect to the Lebesgue measure, of the corresponding distribution.

4.1. *Remarks on the characteristic function of the mean of a gamma process.* Endow $\mathbb{F}$ with the topology of weak convergence, and let $\alpha$ be an element of $\mathbb{F}_a$, with $a > 0$. Denote the Borel $\sigma$-field on $\mathbb{F}$ by $\mathscr{F}$, and call *functional gamma distribution* with parameter $\alpha$ the probability measure $\mathscr{G}_\alpha$ on $(\mathbb{F}, \mathscr{F})$ defined as follows. Set $\tilde{\gamma}$ for the identity map on $\mathbb{F}$ and say that $\tilde{\gamma}$ has the functional gamma distribution with parameter $\alpha$, $\mathscr{G}_\alpha$, if for every finite and measurable partition $\{B_1, \ldots, B_k\}$ of $\mathbb{R}$, the random variables $\tilde{\gamma}(B_1), \ldots, \tilde{\gamma}(B_k)$ are independent, with gamma distribution such that $E(\tilde{\gamma}(B_j)) = \mathrm{Var}(\tilde{\gamma}(B_j)) = \alpha(B_j)$, for each $j = 1, \ldots, k$. It is well known that $\tilde{\gamma}(\cdot)/\tilde{\gamma}(\mathbb{R})$ is a random probability measure with p.d. $\mathscr{D}_\alpha$. Moreover, one easily obtains the following representation for the characteristic function of $\tilde{\gamma}(I) = \int x \tilde{\gamma}(dx)$:

$$(4.1) \qquad \mathscr{G}_\alpha(e^{it\tilde{\gamma}(I)}) = \exp\{-\zeta(-it; \alpha, I)\}, \qquad t \in \mathbb{R}.$$

Thus, according to Section 3.1 in Bondesson ([1992](#)), the p.d. of the random mean $\tilde{\gamma}(I)$ is an extended form of generalized gamma convolution, with *Thorin measure* $\alpha^*$, where $\alpha^* = \alpha \circ J^{-1}$ and $J(x) = 1/x$ for any $x$ in $\mathbb{R} \setminus \{0\}$.

On the other hand, (4.1) can be invoked to establish a representation of the p.d. of $\tilde{\gamma}(I)$ in terms of the p.d. $\mu_\alpha$ of $\tilde{\varphi}(I)$.

PROPOSITION 1. *Suppose that $\alpha$ is an element in $\mathbb{F}_a$, for some $a > 0$, and $0 \le \alpha\{0\} < a$. Then, the p.d. of $\tilde{\gamma}(I)$ is absolutely continuous with respect to the Lebesgue measure and*

$$x \mapsto \mathbb{I}_{(-\infty, 0]}(x) \frac{|x|^{a-1}}{\Gamma(a)} \int_{(-\infty, 0)} |y|^{-a} e^{-x/y} \mu_\alpha(dy)$$

$$+ \mathbb{I}_{(0, +\infty)}(x) \frac{x^{a-1}}{\Gamma(a)} \int_{(0, +\infty)} y^{-a} e^{-x/y} \mu_\alpha(dy) =: q(x)$$

*is a density function of such a p.d.*

PROOF.   Theorem 1(iii) yields

$$\exp\{-\zeta(-it;\alpha,I)\} = \int_{\mathbb{R}}(1-itx)^{-a}\mu_\alpha(dx)$$

$$= \int_{\mathbb{R}}\frac{1}{\Gamma(a)}\int_0^{+\infty}z^{a-1}e^{-z(1-itx)}\,dz\,\mu_\alpha(dx)$$

$$= \int_{-\infty}^0 e^{ity}\frac{|y|^{a-1}}{\Gamma(a)}\int_{(-\infty,0)}|x|^{-a}e^{-y/x}\mu_\alpha(dx)\,dy$$

$$+ \int_0^{+\infty}e^{ity}\frac{y^{a-1}}{\Gamma(a)}\int_{(0,+\infty)}x^{-a}e^{-y/x}\mu_\alpha(dx)\,dy,$$

where the last equality follows from the application of the change-of-variable formula, with $y = y(z) = zx$, and of the Fubini theorem on iterated integrals. It follows that $\exp\{-\zeta(-it;\alpha,I)\} = \mathscr{G}_\alpha(e^{it\tilde{\gamma}(I)})$ is the Fourier transform of the probability density function $q(\cdot)$.   $\square$

Absolute continuity extends to $\mu_\alpha$. We propose a new and simpler proof of this statement. See Regazzini, Guglielmi and Di Nunno (2002) for a different line of reasoning.

PROPOSITION 2.   *Suppose that $\alpha$ is in $\mathbb{F}_a$, for some $a > 0$, with $0 \le \alpha\{x\} < a$ for every $x \in \mathbb{R}$. Then, the p.d. of $\tilde{\varphi}(I)$ under $\mathscr{D}_\alpha$ is absolutely continuous with respect to the Lebesgue measure on $\mathbb{R}$.*

PROOF.   Suppose that $\mu_\alpha$ has a singular part $\mu_{\alpha,s}$ such that $\mu_{\alpha,s}[c,c+\Delta] > 0$ with $\Delta > 0$ and $c$ in $\mathbb{R}$. Since $\alpha$ is nondegenerate, the p.d. function $G_\alpha^c$ of $\tilde{\gamma}(I-c)$ is absolutely continuous and, by virtue of Proposition 1, $G_\alpha^c(\Delta) - G_\alpha^c(0) = (\Gamma(a))^{-1}\int_0^\Delta x^{a-1}\int_{(0,+\infty)}y^{-a}e^{-x/y}\mu_\alpha^c(dy)\,dx$, $\mu_\alpha^c$ being the p.d. of $\tilde{\varphi}(I-c)$. Without real loss of generality, it suffices to consider the case $c = 0$. Then, for any finite class of disjoint intervals $(\rho_1,\rho_1'),\ldots,(\rho_m,\rho_m')$ with $\rho_1 \ge 0$, $\rho_2 > \rho_1',\ldots,\rho_m > \rho_{m-1}'$, and for any $\delta > 0$, one has

$$\sum_{j=1}^m[G_\alpha((1+\delta)\rho_j') - G_\alpha(\rho_j)]$$

$$= \sum_{j=1}^m\int_{\rho_j}^{\rho_j'(1+\delta)}\frac{x^{a-1}}{\Gamma(a)}\int_{(0,+\infty)}\sigma^{-a}e^{-x/\sigma}\mu_\alpha(d\sigma)\,dx$$

$$\ge \sum_{j=1}^m\int_{\rho_j}^{\rho_j'(1+\delta)}\frac{x^{a-1}}{\Gamma(a)}\int_{[\rho_j\vee x/(1+\delta),\rho_j'\wedge x]}\sigma^{-a}e^{-x/\sigma}\mu_\alpha(d\sigma)\,dx$$

$$= \sum_{j=1}^m\int_{\rho_j}^{\rho_j'}\int_1^{1+\delta}\frac{1}{\Gamma(a)}\xi^{a-1}e^{-\xi}\,d\xi\,\mu_\alpha(d\sigma)$$

$$= \int_1^{1+\delta} \frac{1}{\Gamma(a)} \xi^{a-1} e^{-\xi} \sum_{j=1}^m \mu_\alpha[\rho_j, \rho_j'] \, d\xi.$$

Because of the existence of the singular part $\mu_{\alpha,s}$ (with $c = 0$), there is $\varepsilon > 0$ such that: With each $\eta > 0$ one can associate a finite class of disjoint intervals, included in $[0, \Delta]$, $(\rho_1, \rho_1'), \ldots, (\rho_m, \rho_m')$ such that $\sum(\rho_j' - \rho_j) < \eta$ and, for some $\varepsilon > 0$, $\sum \mu_\alpha[\rho_j, \rho_j'] \geq \varepsilon$. On the other hand, from the absolute continuity of $G_\alpha$, there is $\eta$ such that $\sum\{G_\alpha((1+\delta)\rho_j') - G_\alpha(\rho_j)\} < \varepsilon \int_1^{1+\delta} (\Gamma(a))^{-1} \xi^{a-1} e^{-\xi} \, d\xi$ for any finite class of disjoint intervals satisfying $\sum(\rho_j' - \rho_j) < \eta$, yielding a contradiction and, thus, completing the proof.   $\square$

4.2. *Lévy–Khintchine representation of* $\exp\{-\zeta\}$. The next proposition extends part of Theorem 3.1.1 in Bondesson (1992). It involves the distribution function $A$ associated with $\alpha$, the function $g$ defined by

$$
\begin{aligned}
(4.2) \quad g(x) = -\frac{|x|}{1+x^2} \Big[ &\mathbb{I}_{[0,+\infty)}(x) \int_{(0,+\infty)} e^{-yx} \, dA(y^{-1}) \\
&+ \mathbb{I}_{(-\infty,0)}(x) \int_{(-\infty,0)} e^{-yx} \, dA(y^{-1}) \Big], \qquad x \in \mathbb{R}
\end{aligned}
$$

and the well-known fact that the Lévy–Khintchine representation of an infinitely divisible characteristic function, determined by the pair $(\gamma, G)$, is expressed by

$$t \mapsto \exp\Big\{ i\gamma t + \int_{\mathbb{R}\setminus\{0\}} \Big( e^{itu} - 1 - \frac{itu}{1+u^2} \Big) \frac{1+u^2}{u^2} \, dG(u) \Big\},$$

where $\gamma \in \mathbb{R}$, and $G$ is nondecreasing, right-continuous with $G(-\infty) = 0$ and $G(+\infty) < +\infty$.

THEOREM 3. *Suppose $\alpha$ satisfies condition* (ii) *in Theorem* 1. *Then the characteristic function of $\tilde{\gamma}(I)$ is infinitely divisible with Lévy–Khintchine representation determined by the pair* $(\gamma, G)$, *where $\gamma = \int_{\mathbb{R}\setminus\{0\}} |x|^{-1} g(x) \, dx$ and $dG(x) = g(x) \, dx$ for every $x \in \mathbb{R}$, $g$ being defined by* (4.2). *In particular, $\mathscr{G}_\alpha(e^{it\tilde{\gamma}(I)}) \equiv 1$ if and only if $\alpha = \delta_0$.*

PROOF. The result is first proved for the case $\int_{\mathbb{R}} |x| \alpha(dx) < +\infty$. By differentiation under the integral sign one gets

$$
\begin{aligned}
\frac{d}{dt} \log \mathscr{G}_\alpha(e^{it\tilde{\gamma}(I)}) &= \int_{\mathbb{R}} \frac{ix}{1-itx} \alpha(dx) \\
&= \int_{\mathbb{R}} ix \int_0^{+\infty} e^{-z(1-itx)} \, dz \, dA(x)
\end{aligned}
$$

$$= -i \int_0^{+\infty} e^{it\xi} \left( \int_{(0,+\infty)} e^{-\xi y} \, dA(y^{-1}) \right) d\xi$$

$$+ i \int_{-\infty}^0 e^{it\xi} \left( \int_{(-\infty,0)} e^{-\xi y} \, dA(y^{-1}) \right) d\xi,$$

where the application of the Fubini–Tonelli theorem is valid in view of the above extra-assumption. Next, by defining $g$ as in (4.2), one obtains $g \geq 0$, $\int_{\mathbb{R}} g(x) \, dx < +\infty$, $\int_{\mathbb{R}} |x|^{-1} g(x) \, dx < +\infty$. Hence, $G(x) := \int_{-\infty}^x g(u) \, du$ is well defined for every $x$ in $\mathbb{R}$, and it turns out to be nondecreasing with $G(-\infty) = 0$ and $G(+\infty) < +\infty$. Moreover, letting

$$\phi^*(t) := \int_{\mathbb{R}} (e^{itx} - 1) \frac{1+x^2}{x^2} g(x) \, dx, \qquad t \in \mathbb{R},$$

then $e^{\phi^*}$ is an infinitely divisible characteristic function with

$$\gamma = \int_{\mathbb{R} \setminus \{0\}} |x|^{-1} g(x) \, dx$$

and $dG(x) = g(x) \, dx$, such that

$$\frac{d}{dt} \phi^*(t) = \frac{d}{dt} \log \mathscr{G}_\alpha(e^{it\tilde{\gamma}(I)}), \qquad t \in \mathbb{R}.$$

The latter can be verified by interchanging the derivative with the integral in the expression of $\phi^*(\cdot)$. It implies that $e^{\phi^*}$ is the characteristic function of $\tilde{\gamma}(I)$ when $\int_{\mathbb{R}} |x| \alpha(dx) < +\infty$.

The extension of this conclusion to any $\alpha$ such that $\int_{\mathbb{R}} \log(1 + |x|) \alpha(dx) < +\infty$ can be carried out through direct calculation of $\phi^*$. $\quad\square$

REMARK 1. A nice application of this theorem is a further proof of the absolute continuity of the p.d. of $\tilde{\gamma}(I)$. In fact, under the hypotheses in Proposition 2, one has $\int_{\mathbb{R} \setminus \{0\}} x^{-2} g(x) \, dx = +\infty$ and this, applying a criterion due to Tucker (1962), suffices to state the absolute continuity of the (infinitely divisible) p.d. of $\tilde{\gamma}(I)$.

5. Extension of the Markov–Krein identity: the Lauricella identity. In Section 3, the so-called Markov–Krein identity has been obtained as an extension of the Lauricella formulae (2.1) and (2.2) with $a = c > 0$. The Lauricella general representations are now employed to prove an identity that holds for all pairs $(a, c)$ of strictly positive numbers. For the remainder of this section,

$$I_1 := \int_0^{(1+)} \exp\{-\zeta(iwt; \alpha, I)\} b(w; c, a - c) \, dw$$

is to be meant as the integral along the contour $OGG'$, as shown in Figure 1, of the function $\exp\{-\zeta(iwt;\alpha,I)\}b(w;c,a-c)$, where

$$b(w;c,a-c) = \frac{\Gamma(a)}{\Gamma(c)\Gamma(a-c)}w^{c-1}(1-w)^{a-c-1}$$

for every $w$ in $\mathbb{C}$ such that $\operatorname{Re}(w) > 0$ and $\operatorname{Im}(w) \neq 0$ if $\operatorname{Re}(w) \in (0,1]$.

In view of the Cauchy integral theorem, this contour can be deformed into the path of integration consisting of: (A) the straight line segment $w = x - i\tau$ with $x$ varying in $(0, 1-\varepsilon)$, (C) the circle $1 + \varepsilon e^{i\theta}$ with $\theta \in (-\pi + \eta, \pi - \eta)$, for a suitable $\eta > 0$ and (A') the straight line segment $w = x + i\tau$. Thus, if $c$ and $(a - c)$ are strictly positive, with $(a - c)$ different from $1, 2, \dots$, one can let $\varepsilon \to 0$ to obtain

$$(5.1) \qquad I_1 = (1 - e^{2\pi i(a-c)})\int_0^1 \exp\{-\zeta(iwt;\alpha,I)\}b(w;c,a-c)\,dw.$$

This is required to state the following extension of point (iii) in Theorem 1, which we shall call *Lauricella's identity*. The symbols $B$ and $\mu_\alpha$ appearing in this identity have the same meaning as in Sections 2 and 3, respectively.

THEOREM 4. *If, for some $a > 0$, $\alpha$ is in $\mathbb{F}_a$, then*

$$(5.2) \qquad \int_{\mathbb{R}}(1 + itx)^{-c}\mu_\alpha(dx) = \int_{[0,1]}\exp\{-\zeta(iut;\alpha,I)\}B(\,du;c,a-c)$$

*if $a \geq c > 0$ and*

$$(5.3) \qquad \begin{aligned} &\int_{\mathbb{R}}(1 + itx)^{-c}\mu_\alpha(dx) \\ &\quad = \frac{\Gamma(c-a+1)\Gamma(a)}{2\pi i\Gamma(c)}\int_0^{(1+)}\exp\{-\zeta(iwt;\alpha,I)\}w^{c-1}(w-1)^{a-c-1}\,dw \end{aligned}$$

*if $c > a > 0$.*

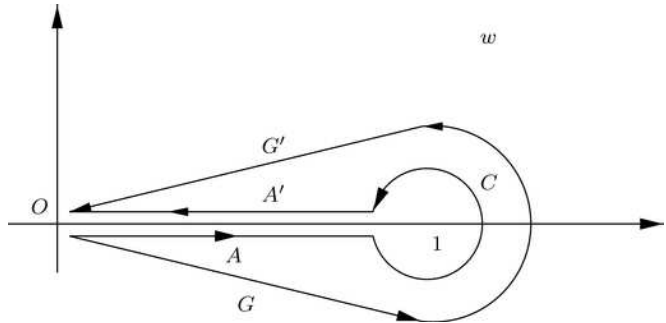

FIG. 1.

PROOF.  First note that, if $0 < c \leq a$, the following equality holds true by virtue of analytic continuation of both sides of (2.1) and (2.2):

$$\int_{\mathbb{P}} \frac{1}{(1 + it\varphi(\xi_n))^c} \mathscr{D}_\alpha(\, d\varphi) = \int_{[0,1]} \exp\{-\zeta(iut; \alpha_n, I)\} B(\, du; c, a - c).$$

Next, argue as in the proof of Theorem 1 and apply the basic limit theorems of integration theory to obtain (5.2).

Finally, suppose that $a > c > 0$ and use what has just been proved and (5.1) to state the equality

$$\int_{\mathbb{R}} \frac{1}{(1 + itx)^c} \mu_\alpha(dx) = \frac{I_1}{1 - e^{2\pi i(a-c)}}$$

when $a - c \neq 1, 2, \ldots$. Now, take Euler's reflection formula $\Gamma(z)\Gamma(1 - z) = \pi(\sin \pi z)^{-1}$ to obtain

$$\int_{\mathbb{R}} \frac{1}{(1 + itx)^c} \mu_\alpha(dx)$$

$$= \frac{\Gamma(c - a + 1)\Gamma(a)}{2\pi i \Gamma(c)} \int_0^{(1+)} w^{c-1}(w - 1)^{a-c-1} \exp\{-\zeta(iwt; \alpha, I)\}\, dw.$$

Notice that the right-hand side of the previous equality, as a function of $c$, is analytic on $\{c \in \mathbb{C} : \mathrm{Re}(c) > 0, \mathrm{Re}(a - c) \notin \{0, 1, 2, \ldots\}\}$. Hence, such an equality can be extended to any $c$ with $c > a$. This completes the proof.  □

The Lauricella identity establishes that the Stieltjes transform (of order $c > 0$) of $\mu_\alpha$ can be viewed as a mixture of the Stieltjes transforms (of order $a$) of p.d.'s of means of Dirichlet processes with parameter $\alpha_w(\cdot) := \alpha\{x : wx \in \cdot\}$. Apropos of this, recall that

$$\frac{\Gamma(c - a + 1)\Gamma(a)}{2\pi i \Gamma(c)} \int_0^{(1+)} w^{c-1}(w - 1)^{a-c-1}\, dw = 1.$$

Compare, for example, 3.1.27 in Slater (1960).

## 6. Applications of the previous results.

6.1. *Exact forms of the p.d. of $\int x\tilde{\varphi}(dx)$.*  Much work has been done on the exact form of the p.d. of $\tilde{\varphi}(I)$. As the following arguments show, Lauricella's identity considerably simplifies the solution of the problem. First, rewrite Theorem 4 with $c = 1$ and $\alpha$ such that $\alpha\{x\} < \alpha(\mathbb{R})$ for every $x$. Denoting a density function for $\mu_\alpha$ by $m_\alpha$ and setting $L(w) :=$

$\int_{\mathbb{R}}(w+x)^{-1}m_\alpha(x)\,dx$, one has

$$L(w) = \begin{cases} w^{-1}\int_{[0,1]}\exp\{-\zeta(uw^{-1};\alpha,I)\}B(\,du;1,a-1), & \text{if } a>1, \\[2mm] w^{-1}\exp\{-\zeta(w^{-1};\alpha,I)\}, & \text{if } a=1, \\[2mm] w^{-1}(1-\exp\{2\pi i(a-1)\})^{-1} \\[2mm] \quad \times \int_0^{(1+)}\exp\{-\zeta(uw^{-1};\alpha,I)\}b(u;1,a-1)\,du, & \text{if } a\in(0,1), \end{cases}$$

for every $w\in\mathbb{C}$ for which $\mathrm{Im}(w)\neq 0$. Thus, the Stieltjes–Perron inversion formula yields

$$\mu_\alpha(x_1,x_2) = \frac{1}{\pi}\lim_{\varepsilon\downarrow 0}\mathrm{Im}\int_{x_1}^{x_2}\frac{1}{(-\lambda-i\varepsilon)}\exp\left\{-\zeta\left(\frac{1}{-\lambda-i\varepsilon};\alpha,I\right)\right\}d\lambda$$

if $a=1$,

$$\begin{aligned} \mu_\alpha(x_1,x_2) & \\ = \frac{a-1}{\pi}\lim_{\varepsilon\downarrow 0}&\mathrm{Im}\int_{x_1}^{x_2}\frac{1}{(-\lambda-i\varepsilon)} \\ & \times\int_0^1\exp\left\{-\zeta\left(\frac{u}{-\lambda-i\varepsilon};\alpha,I\right)\right\}(1-u)^{a-2}\,du\,d\lambda \end{aligned}$$

if $a>1$, and

$$\begin{aligned} \mu_\alpha(x_1,x_2) & \\ = \frac{1}{\pi}\lim_{\varepsilon\downarrow 0}&\mathrm{Im}\int_{x_1}^{x_2}\frac{1}{(1-e^{2\pi i(a-1)})(-\lambda-i\varepsilon)} \\ & \times\int_0^{(1+)}\exp\left\{-\zeta\left(\frac{w}{-\lambda-i\varepsilon};\alpha,I\right)\right\}b(w;1,a-1)\,dw\,d\lambda \end{aligned}$$

if $a\in(0,1)$. In particular, when $a=1$, one obtains

$$m_\alpha(\xi) = \frac{1}{\pi}\lim_{\varepsilon\downarrow 0}\mathrm{Im}\frac{1}{(-\xi-i\varepsilon)}\exp\left\{-\zeta\left(\frac{1}{-\xi-i\varepsilon};\alpha,I\right)\right\};$$

if $a>1$ and the saltus of $A$ at each discontinuity point is strictly smaller than 1, one has

$$m_\alpha(\xi) = \frac{1}{\pi}\lim_{\varepsilon\downarrow 0}\mathrm{Im}\frac{1}{(-\xi-i\varepsilon)}\int_{[0,1]}\exp\left\{-\zeta\left(\frac{u}{-\xi-i\varepsilon};\alpha,I\right)\right\}B(\,du;1,a-1).$$

At this stage we can achieve the same results as in Proposition 9(ii) and (iii) of Regazzini, Guglielmi and Di Nunno (2002), by a simple application of the change-of-variable formula.

EXAMPLE 1. (a) Let $\alpha(dx) = C_{\theta,\sigma^2}(dx)$ be the Cauchy density function with parameters $\theta \in \mathbb{R}$ and $\sigma > 0$, that is,

$$C_{\theta,\sigma^2}(dx) := \frac{\sigma}{\pi(1 + \sigma^2(x - \theta)^2)}\, dx, \qquad x \in \mathbb{R}.$$

It is also convenient to identify $C_{\theta,+\infty}$ with the degenerate distribution $\delta_\theta$. By resorting to the expression of $m_\alpha$ given above when $a = 1$, it is easily verified that, for any $x$ in $\mathbb{R}$,

$$
\begin{aligned}
m_\alpha(x) &= \frac{1}{\pi}\exp\left\{-\int_{\mathbb{R}\setminus\{x\}}\log|\rho - x|C_{\theta,\sigma}(\,d\rho)\right\}\sin\left(\frac{\pi}{2} + \arctan(\sigma(x - \theta))\right) \\
&= \frac{\sigma}{\pi(1 + \sigma^2(x - \theta)^2)}
\end{aligned}
$$

if $\sigma \in (0 + \infty)$, whereas $\mu_\alpha$ is degenerate at $\theta$ if $\sigma = +\infty$. In other terms, if $\alpha$ coincides with $C_{\theta,\sigma^2}$ for some $\theta \in \mathbb{R}$ and $\sigma \in (0, +\infty]$, then $\mu_\alpha = \alpha$, a well-known result stated by Yamato ([1984](#)). In Section 6.2 the converse will be proved.

(b) If $\alpha(dx) = a\mathbb{I}_{(0,1)}(x)\, dx$ with $a > 1$, then, for any $\xi \in (0,1)$,

$$
\begin{aligned}
m_\alpha(\xi) &= \frac{(a-1)(1-\xi)^{a-1}e^{a\xi}}{\pi} \\
&\quad \times \int_0^1 \frac{e^{a(1-\xi)u}u^{a-2}\sin(a\pi(1-\xi)(1-u))}{[(1-\xi)u + \xi]^{a\xi + a(1-\xi)u}[(1-\xi)(1-u)]^{a(1-\xi)(1-u)}}\, du.
\end{aligned}
$$

(c) If $\alpha(dx) = a(\sigma\sqrt{2\pi})^{-1}\exp\{-\frac{1}{2\sigma^2}(x-\theta)^2\}\, dx$, for any $\sigma > 0$, $\theta \in \mathbb{R}$ and $a > 1$, then application of 2.6.22.1 in Prudnikov, Brychkov and Marichev ([1986](#)) yields

$$
\begin{aligned}
m_\alpha(\xi) &= \frac{(a-1)\xi^{a-1}}{\pi} \\
&\quad \times \int_0^1 (1-u)^{a-2}\sin\left[a\pi\left\{1 - \Phi\left(\frac{\xi - u\theta}{u\sigma}\right)\right\}\right] \\
&\quad \times \exp\left\{-\frac{a}{\sigma\sqrt{2\pi}}\exp\left\{-\frac{(\xi - u\theta)^2}{2u^2\sigma^2}\right\}\right. \\
&\quad \left.\times \frac{\partial}{\partial\nu}\left((u\sigma)^\nu 2^{-(\nu/2)}\Gamma(\nu)\Psi\left(\frac{\nu}{2}, \frac{1}{2}; \frac{(\xi - u\theta)^2}{2\sigma^2}\right)\right)\Big|_{\nu=1}\right\}\, du,
\end{aligned}
$$

for any $\xi > 0$, where $\Phi$ denotes the distribution function of a Gaussian random variable with zero mean and variance equal to 1 and $\Psi$ is the Tricomi confluent hypergeometric function. On the other hand, letting $\xi < 0$, one has

$$m_\alpha(\xi) = -\frac{(a-1)|\xi|^{a-1}}{\pi}$$

$$\times \int_0^1 (1-u)^{a-2} \sin\left[a\pi\Phi\left(\frac{\xi-u\theta}{u\sigma}\right)\right]$$

$$\times \exp\left\{-\frac{a}{\sigma\sqrt{2\pi}}\exp\left\{-\frac{(\xi-u\theta)^2}{2u^2\sigma^2}\right\}\right.$$

$$\left.\times \frac{\partial}{\partial\nu}\left((u\sigma)^\nu 2^{-(\nu/2)}\Psi\left(\frac{\nu}{2},\frac{1}{2};\frac{(\xi-u\theta)^2}{2\sigma^2}\right)\right)\Big|_{\nu=1}\right\}\,du.$$

6.2. *Characterization of the Cauchy distribution.* We prove the characterization of the Cauchy distribution already mentioned in Example 1(a). The proof is based essentially on the result in Theorem 4. As in Theorem 2, $\beta_a$ is the bijection of $\mathbb{M}_a$ on $\mathbb{F}_a$.

THEOREM 5. *The class of all fixed points of the bijection $\beta_1$ is $\{C_{\theta,\sigma^2} : \theta \in \mathbb{R}, \sigma^2 \in (0, +\infty]\}$.*

In other terms, if $\alpha \in \mathbb{F}_1$, then $\mu_\alpha = \alpha$ if and only if $\alpha$ is Cauchy or degenerate. In Section 4 of Cifarelli and Regazzini (1990) another characterization is given under stronger hypotheses, that is: A nondegenerate probability is a fixed point of $\beta_a$ *for every $a$* if and only it is Cauchy.

PROOF OF THEOREM 5. In the light of Example 1(a), it suffices to prove that the condition is necessary. Suppose that $\mu_\alpha = \alpha \in \mathbb{F}_1$. Then $\mu_{n\alpha} = \alpha$ for every $n \in \mathbb{N}$. Indeed, by the Markov–Krein identity,

$$(6.1) \qquad \int_\mathbb{R} \frac{1}{s+iy}\alpha(dy) = \exp\left\{-\int_\mathbb{R}\log(s+iy)\alpha(dy)\right\}$$

holds for every $s \in \mathbb{R}\setminus\{0\}$ by virtue of the above hypothesis. Differentiate both sides to obtain

$$\int_\mathbb{R}\frac{1}{(s+iy)^2}\alpha(dy) = \exp\left\{-\int_\mathbb{R}\log(s+iy)\alpha(dy)\right\}\int_\mathbb{R}\frac{1}{s+iy}\alpha(dy)$$

$$= \exp\left\{-2\int_\mathbb{R}\log(s+iy)\alpha(dy)\right\} \qquad [\text{by (6.1)}].$$

Hence, the result holds true for $n = 2$. Assume that it is valid for all $m < n$. Then

$$\int_\mathbb{R}\frac{1}{(s+iy)^{n-1}}\alpha(dy) = \exp\left\{-(n-1)\int_\mathbb{R}\log(s+iy)\alpha(dy)\right\}$$

and, by differentiation,

$$\int_{\mathbb{R}} \frac{1}{(s+iy)^n} \alpha(dy) = \exp\left\{-(n-1)\int_{\mathbb{R}} \log(s+iy)\alpha(dy)\right\} \int_{\mathbb{R}} \frac{1}{s+iy}\alpha(dy)$$

$$= \exp\left\{-n\int_{\mathbb{R}} \log(s+iy)\alpha(dy)\right\} \qquad \text{[by (6.1)]}.$$

By induction one can conclude that $\mu_{n\alpha} = \alpha$ for every $n \geq 1$. Next, this statement and Theorem 4, with $a=2$ and $c=1$, give

$$\frac{1}{n-1}\int_{\mathbb{R}} \frac{1}{1+itx}\mu_{n\alpha}(dx) = \int_0^1 \exp\{-\zeta(iut;n\alpha,I)\}(1-u)^{n-2}\,du$$

where, by hypothesis,

$$\int_{\mathbb{R}} \frac{1}{1+itx}\mu_{n\alpha}(dx) = \int_{\mathbb{R}} \frac{1}{1+itx}\mu_{\alpha}(dx) = \exp\{-\zeta(it;\alpha,I)\}.$$

Hence,

$$\frac{t}{n-1}\exp\{-\zeta(it;\alpha,I)\} = \int_0^t \exp\{-\zeta(iu;n\alpha,I)\}\left(1-\frac{u}{t}\right)^{n-2}\,du.$$

In particular, letting $\zeta(it;\alpha,I) = \rho(t)$ and $n=2$, one has

$$te^{-\rho(t)} = \int_0^t e^{-2\rho(x)}\,dx, \qquad t \in \mathbb{R}.$$

Differentiation of both sides gives $e^{-\rho(t)} - t\rho'(t)e^{-\rho(t)} = e^{-2\rho(t)}$, that is,

$$\rho(t) = \log(1+wt)$$

for some $w \in \mathbb{C}$. Hence, $\exp\{-\int_{\mathbb{R}} \log(s+ix)\alpha(dx)\} = (s+w)^{-1}$ holds true for every $s > 0$, and by (6.1) one gets

$$(s+w)^{-1} = \int_{\mathbb{R}} (s+iy)^{-1}\alpha(dy) = \int_0^{+\infty} e^{-sx}\left(\int_{\mathbb{R}} e^{-i\xi y}\alpha(dy)\right)\,d\xi.$$

The left-hand side is the moment generating function of $e^{-wx}\mathbb{I}_{(0,+\infty)}(x)$ for every $s > -\text{Re}(w)$, implying $e^{-wx} = \int_{\mathbb{R}} e^{-ix y}\alpha(dy)$, for every $x > 0$. Analogously, $e^{wx} = \int_{\mathbb{R}} e^{-ixy}\alpha(dy)$ is valid for any $x < 0$. Therefore,

$$\alpha(dy) = C_{\theta,\sigma^2}(dy), \qquad y \in \mathbb{R},$$

where $\sigma = 1/\text{Im}(w)$, $\text{Im}(w) \geq 0$ and $\theta = \text{Re}(w)$.  $\square$

REMARK 2. Richard Olshen has drawn the attention of the authors to the following problem. Let $(X_n)_{n\geq 1}$ be a sequence of exchangeable real-valued random variables, whose p.d. $P$ is supposed to be characterized by

$$(6.2) \qquad P(A_1 \times \cdots \times A_n \times \mathbb{R}^{\infty}) = \int_{\mathbb{P}} p(A_1)\cdots p(A_n)\mathscr{D}_{\alpha}(\,dp),$$

$n \in \mathbb{N}$, $A_k \in \mathscr{B}(\mathbb{R})$ for $k = 1, \ldots, n$, $\alpha$ being any element of $\mathbb{F}_1$. Then, by de Finetti's representation theorem and the strong law of large numbers for i.i.d. random variables, one has

$$P\left(\left|\sum_{k=1}^{n} X_k/n - \tilde{\varphi}(I)\right| \to 0 \,\Big|\, \tilde{\varphi}\right) = 1 \qquad \text{a.s.-}P,$$

which entails

$$\frac{1}{n}\sum_{k=1}^{n} X_k \to \int_{\mathbb{R}} x\tilde{\varphi}(dx) \qquad \text{a.s.-}P.$$

Thus, the limiting p.d. of $\sum_{k=1}^{n} X_k/n$ coincides with the law of $\tilde{\varphi}(I)$, that is, with $\mu_\alpha$. Moreover, recall that each $X_k$ has p.d. $\alpha$. A natural question is whether the above limiting p.d. may coincide with $\alpha$. An answer can be deduced from Theorem 5, that is: *If $(X_n)_{n \geq 1}$ is a sequence of exchangeable random variables with a law characterized by* (6.2), *then the limiting distribution of the empirical mean $\sum_{k=1}^{n} X_k/n$ is $\alpha$ if and only if $\alpha$ is an element of $\{C_{\theta,\sigma^2} : \theta \in \mathbb{R}, \sigma^2 \in (0, +\infty]\}$.*

**7. Characteristic function of $\int x\tilde{\varphi}(dx)$.** Let $b_1, \ldots, b_k$ be strictly positive numbers such that $|\mathbf{b}| < a$, and let $x_1, \ldots, x_n$ be arbitrary real numbers. Erdélyi (1937) defines a confluent form of the fourth Lauricella function, ${}_n\Phi$, through the following limiting process:

$${}_n\Phi(b_1, \ldots, b_n; a; itx_1, \ldots, itx_n) = \lim_{\varepsilon \downarrow 0} F_D(\varepsilon^{-1}; b_1, \ldots, b_n; a; i\varepsilon tx_1, \ldots, i\varepsilon tx_n),$$

where $t$ is any real number. From (8.5) in Erdélyi (1937) and the definition of $\mathscr{D}_\alpha$, we get

$${}_n\Phi(b_1, \ldots, b_n; a; itx_1, \ldots, itx_n)$$

$$= \frac{\Gamma(a)}{\Gamma(b_1)\cdots\Gamma(b_n)\Gamma(a - |\mathbf{b}|)}$$

$$\times \int_{T_n} u_1^{b_1-1}\cdots u_n^{b_n-1}(1 - |\mathbf{u}|)^{a-|\mathbf{b}|-1}\exp\left\{it\sum_{k=1}^{n} x_k u_k\right\} du_1\cdots du_n$$

$$= \mathscr{D}_\alpha\left(\exp\left\{it\int x\tilde{\varphi}(dx)\right\}\right),$$

where $\tilde{\varphi}$ is a Dirichlet process with parameter $\alpha$ defined by $\alpha\{x_k\} = b_k$, $k = 1, \ldots, n$, and $\alpha\{0\} = a - |\mathbf{b}|$. Hence, ${}_n\Phi$ can be viewed as the characteristic function of the mean of a Dirichlet process with a parameter $\alpha$ supported by a finite set.

The function $_n\Phi$ admits an interesting representation, as a single contour integral, which Erdélyi achieved by resorting to the following argument. Take the Laplace transform of

$$\varphi^*(\sigma) := (\Gamma(a))^{-1}\sigma^{a-1}{}_n\Phi(b_1,\ldots,b_n;a;it\sigma x_1,\ldots,it\sigma x_n),$$

namely,

$$\int_0^{+\infty}\varphi^*(\sigma)e^{-w\sigma}\,d\sigma$$

$$= \frac{1}{w^a}\mathscr{D}_\alpha\left[\frac{w^a}{\Gamma(a)}\int_0^{+\infty}\sigma^{a-1}\exp\left\{-w\sigma+it\sigma\int x\tilde{\varphi}(dx)\right\}d\sigma\right]$$

$$\qquad\qquad\qquad\qquad\qquad\qquad [\operatorname{Re}(w) > 0]$$

$$= \int_{\mathbb{P}}\frac{1}{(w-it\int x\varphi(dx))^a}\mathscr{D}_\alpha(\,d\varphi)$$

$$= \frac{1}{w^a}\exp\left\{-\int_{\mathbb{R}}\log\left(1-\frac{itx}{w}\right)\alpha(dx)\right\}\qquad[\text{by Theorem 1(iii)}]$$

$$= \frac{1}{w^a}\prod_{k=1}^n\left\{1-\frac{itx_k}{w}\right\}^{-b_k}.$$

The well-known complex inversion formula for the Laplace transform [see, e.g., Henrici (1991), page 278] yields

$$\frac{x^{a-1}}{\Gamma(a)}\mathscr{D}_\alpha\left(e^{it\tilde{\varphi}(I)}\right) = \varphi^*(x)$$

$$= \lim_{R\to+\infty}\frac{1}{2\pi i}\int_{\gamma-iR}^{\gamma+iR}e^{wx}\frac{1}{w^a}\prod_{k=1}^n\left\{1-\frac{itx_k}{w}\right\}^{-b_k}dw\qquad\gamma > 0,$$

which, letting $x = 1$, gives

$$\mathscr{D}_\alpha\left(e^{it\tilde{\varphi}(I)}\right) = \frac{\Gamma(a)}{2\pi i}\lim_{R\to+\infty}\int_{\gamma-iR}^{\gamma+iR}e^w\frac{1}{w^a}\prod_{k=1}^n\left\{1-\frac{itx_k}{w}\right\}^{-b_k}dw$$

$$= \frac{\Gamma(a)}{2\pi i}\text{PV}\int_{\gamma-i\infty}^{\gamma+i\infty}e^w\frac{1}{w^a}\exp\left\{-\zeta\left(-\frac{it}{w};\alpha,I\right)\right\}dw,$$

PV$\int$ denoting *principal value integral*. Moreover, from the Cauchy integral theorem, the path of integration can be deformed into any contour which consists of a simple loop, $\mathfrak{C}_t$, beginning and ending at $-\infty$, and encircling all the finite singularities of the integrand, that is: $0, itx_1,\ldots,itx_n$. See Figure 2.

The above discussion leads to a first representation—under some restrictions—for the characteristic function of $\tilde{\varphi}(I)$.

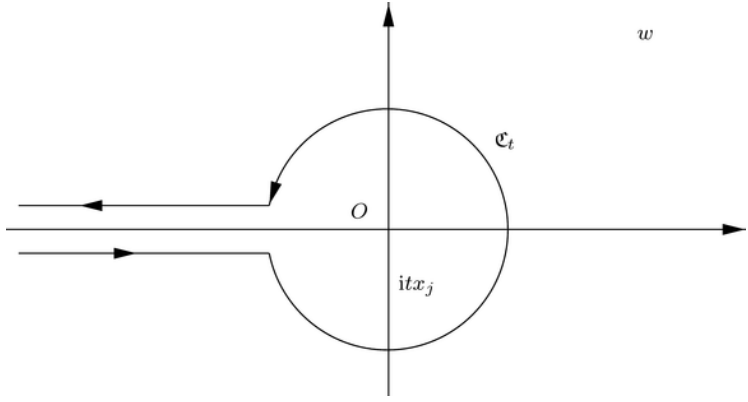

Fig. 2.

THEOREM 6.  *Let $\alpha$ be a measure with a bounded support $S \subset \mathbb{R}$, and, for any $t$ in $\mathbb{R}$, let $S_t$ be the closure of the convex hull of $\{itx : x \in S\}$. Then, with $a = \alpha(S) > 0$,*

$$\int_{\mathbb{P}} e^{it\varphi(I)} \mathscr{D}_\alpha(\,d\varphi) = \frac{\Gamma(a)}{2\pi i} \mathrm{PV} \int_{\gamma-i\infty}^{\gamma+i\infty} e^w \frac{1}{w^a} \exp\Big\{-\zeta\Big(-\frac{it}{w}; \alpha, I\Big)\Big\}\,dw$$

*holds for every $\gamma > 0$. Moreover, if $\mathfrak{C}_t$ is any contour consisting of a simple loop, beginning and ending at $-\infty$, and encircling $S_t$ in such a way that $S_t \cap \mathfrak{C}_t = \varnothing$, then*

$$\int_{\mathbb{P}} e^{it\varphi(I)} \mathscr{D}_\alpha(\,d\varphi) = \frac{\Gamma(a)}{2\pi i} \int_{\mathfrak{C}_t} \frac{e^w}{w^a} \exp\Big\{-\zeta\Big(-\frac{it}{w}; \alpha, I\Big)\Big\}\,dw.$$

PROOF.  In view of the previous discussion, the statements are true if $S$ is finite. More in general, if $S$ is bounded, resort to the approximating functions $\xi_n$, used in the proof of Theorem 1, to write

$$\int_{\mathbb{P}} e^{it\varphi(I)} \mathscr{D}_\alpha(\,d\varphi) = \lim_{n\to+\infty} \int_{\mathbb{P}} e^{it\varphi(\xi_n)} \mathscr{D}_\alpha(\,d\varphi)$$
$$= \frac{\Gamma(a)}{2\pi i} \lim_{n\to+\infty} \int_{\mathfrak{C}_t} \frac{e^w}{w^a} \exp\Big\{-\zeta\Big(-\frac{it}{w}; \alpha_n, I\Big)\Big\}\,dw.$$

Next, interchanging limit and integral gives

$$\lim_{n\to+\infty} \int_{\mathfrak{C}_t} \frac{e^w}{w^a} \exp\Big\{-\zeta\Big(-\frac{it}{w}; \alpha_n, I\Big)\Big\}\,dw = \int_{\mathfrak{C}_t} \frac{e^w}{w^a} \exp\Big\{-\zeta\Big(-\frac{it}{w}; \alpha, I\Big)\Big\}\,dw.$$

Finally, by the Cauchy integral theorem, one has

$$\int_{\mathfrak{C}_t} \frac{e^w}{w^a} \exp\Big\{-\zeta\Big(-\frac{it}{w}; \alpha, I\Big)\Big\}\,dw = \mathrm{PV} \int_{\gamma-i\infty}^{\gamma+i\infty} \frac{e^w}{w^a} \exp\Big\{-\zeta\Big(-\frac{it}{w}; \alpha, I\Big)\Big\}\,dw$$

for every $\gamma > 0$.  $\square$

To deal with general parameters $\alpha$, assume that condition (ii) in Theorem 1 is satisfied and, for any $k$ in $\mathbb{N}$, define $\alpha^{(k)}$ as

$$\alpha^{(k)}(\cdot) = \delta_{-k}(\cdot)\alpha(-\infty, -k] + \alpha(\cdot \cap (-k, k)) + \delta_k(\cdot)\alpha[k, +\infty).$$

It can be proved that the distribution of $\int x \tilde{\varphi}_k(dx)$, where $\tilde{\varphi}_k$ is the Dirichlet process with parameter $\alpha^{(k)}$, converges weakly to the mean of the Dirichlet process with parameter $\alpha$, as $k$ tends to $+\infty$. Hence, by the continuity theorem,

$$(7.1) \qquad \int_{\mathbb{P}} e^{it\varphi(I)} \mathscr{D}_{\alpha^{(k)}}(d\varphi) \to \int_{\mathbb{P}} e^{it\varphi(I)} \mathscr{D}_\alpha(d\varphi), \qquad t \in \mathbb{R}, k \to +\infty.$$

From Theorem 6 it follows

$$\int_{\mathbb{P}} e^{it\varphi(I)} \mathscr{D}_{\alpha^{(k)}}(d\varphi) = \frac{\Gamma(a)}{2\pi i} \mathrm{PV} \int_{\gamma - i\infty}^{\gamma + i\infty} \frac{e^w}{w^a} \exp\left\{-\zeta\left(-\frac{it}{w}; \alpha^{(k)}, I\right)\right\} dw,$$
(7.2)

which proves

THEOREM 7. *Let $\alpha$ be a measure in $\mathbb{F}_a$, for some $a > 0$. Then*

$$(7.3) \quad \begin{aligned} &\int_{\mathbb{P}} e^{it\varphi(I)} \mathscr{D}_\alpha(d\varphi) \\ &= \lim_{k \to +\infty} \frac{\Gamma(a)}{2\pi i} \mathrm{PV} \int_{\gamma - i\infty}^{\gamma + i\infty} \frac{e^w}{w^a} \exp\left\{-\zeta\left(-\frac{it}{w}; \alpha^{(k)}, I\right)\right\} dw \end{aligned}$$

*holds for every $t$ in $\mathbb{R}$ and for any strictly positive $\gamma$.*

REMARK 3. Principal value integral in the right-hand side of (7.3) reduces to Lebesgue integral if $a > 1$. This is the case when, in the presence of exchangeable observations, one is interested in the Fourier transform of any conditional p.d. of $\tilde{\varphi}(I)$, given those observations.

As a first application of this result, consider the problem of characterizing symmetric p.d.'s in $\mathbb{M}_a$ for some fixed $a > 0$. It is well known that the p.d. of $\tilde{\varphi}(I)$, under $\mathscr{D}_\alpha$, is symmetric if $\alpha$ is symmetric. See, for example, Regazzini, Guglielmi and Di Nunno (2002). Here, in addition to a simple proof of this fact, we prove the converse.

THEOREM 8. *Let $\alpha$ be an element in $\mathbb{F}_a$. Then the distribution $\mu_\alpha$ of the mean of a Dirichlet process with parameter $\alpha$ is symmetric if and only if $\alpha$ is symmetric.*

Proof.   Suppose that $\alpha$ in $\mathbb{F}_a$ is symmetric. Then, for any $k \in \mathbb{N}$,

$$\frac{\Gamma(a)}{2\pi i}\text{PV}\int_{\gamma-i\infty}^{\gamma+i\infty}\frac{e^w}{w^a}\exp\left\{-\int_{\mathbb{R}}\log\left(1-\frac{itx}{w}\right)dA^{(k)}(x)\right\}dw$$

$$=\frac{\Gamma(a)}{2\pi i}\text{PV}\int_{\gamma-i\infty}^{\gamma+i\infty}\frac{e^w}{w^a}\exp\left\{-\int_{\mathbb{R}}\log\left(1-\frac{itx}{w}\right)d(a-A^{(k)}(-x))\right\}dw,$$

which, by the change of variable $y=-x$, can be shown to be equal to

$$\frac{\Gamma(a)}{2\pi i}\text{PV}\int_{\gamma-i\infty}^{\gamma+i\infty}\frac{e^w}{w^a}\exp\left\{-\int_{\mathbb{R}}\log\left(1+\frac{ity}{w}\right)dA^{(k)}(y)\right\}dw,$$

thus implying, by virtue of Theorem 7, that the characteristic function of $\tilde{\varphi}(I)$ is real-valued.
Conversely, if $\mu_\alpha$ in $\mathbb{M}_a$ is symmetric, one has

$$\lim_{k\to+\infty}\frac{\Gamma(a)}{2\pi i}\text{PV}\int_{\gamma-i\infty}^{\gamma+i\infty}\frac{e^w}{w^a}\exp\left\{-\int_{\mathbb{R}}\log\left(1-\frac{itx}{w}\right)d(a-A^{(k)}(-x))\right\}dw$$

$$=\lim_{k\to+\infty}\frac{\Gamma(a)}{2\pi i}\text{PV}\int_{\gamma-i\infty}^{\gamma+i\infty}\frac{e^w}{w^a}\exp\left\{-\int_{\mathbb{R}}\log\left(1+\frac{itx}{w}\right)dA^{(k)}(x)\right\}dw$$

$$\text{(by change of variable)}$$

$$=\lim_{k\to+\infty}\frac{\Gamma(a)}{2\pi i}\text{PV}\int_{\gamma-i\infty}^{\gamma+i\infty}\frac{e^w}{w^a}\exp\left\{-\int_{\mathbb{R}}\log\left(1-\frac{itx}{w}\right)dA^{(k)}(x)\right\}dw$$

$$\text{(by symmetry of }\mu_\alpha).$$

The proof is concluded by resorting to Theorem 2 (uniqueness).   $\square$

Another application concerns the problem of determining the p.d. of $\tilde{\varphi}(f)=\int f\,d\tilde{\varphi}$, when $\tilde{\varphi}$ has the functional Dirichlet distribution $\mathscr{D}_\alpha$, and $f$ is any measurable real-valued function. Here, $\mathscr{D}_\alpha$ can be understood as the functional Dirichlet distribution of a random probability measure $\tilde{m}$ on a probability space $(\Omega,\mathscr{F})$; that is, $\Omega$ and $\mathscr{F}$ need not coincide with $\mathbb{R}$ and $\mathscr{B}(\mathbb{R})$, respectively. In this general framework, $\alpha$ must be a finite measure on $(\Omega,\mathscr{F})$ satisfying $a=\alpha(\Omega)>0$. If $f$ is bounded, then the discretization process used to prove Theorem 6 can be applied to give

$$\mathscr{D}_\alpha(e^{it\tilde{m}(f)})=\frac{\Gamma(a)}{2\pi i}\text{PV}\int_{\gamma-i\infty}^{\gamma+i\infty}\frac{e^w}{w^a}\exp\left\{-\int_{\mathbb{R}}\log\left(1-\frac{itx}{w}\right)\alpha\circ f^{-1}(dx)\right\}dw,$$

allowing an obvious generalization of Theorem 8.

Theorem 9.   *Let $\tilde{m}$ be a random probability measure on $(\Omega,\mathscr{F})$ with distribution $\mathscr{D}_\alpha$, and let $f$ be a measurable function from $\Omega$ to $\mathbb{R}$ such that $\int_\Omega\log(1+|f|)\,d\alpha<+\infty$. Then, the p.d. of $\tilde{m}(f)$ is the same as the p.d. of $\tilde{\varphi}(I)$, where $\tilde{\varphi}$ is a random probability measure on $(\mathbb{R},\mathscr{B}(\mathbb{R}))$ with p.d. $\mathscr{D}_{\alpha\circ f^{-1}}$.*

**8. Vector of means of a single Dirichlet process.** Let $(\Omega, \mathscr{F}, P)$ and $\mathscr{D}_\alpha$ be those defined in the final part of the previous section, and let $f_1, \ldots, f_d$ be measurable functions from $\Omega$ to $\mathbb{R}$ satisfying

$$(8.1) \qquad \int_\Omega \log(1 + |f_k|) \, d\alpha < +\infty, \qquad k = 1, \ldots, d.$$

For any $\mathbf{t} = (t_1, \ldots, t_d)$ in $\mathbb{R}^d$, define $\alpha_{\langle \mathbf{t}, \mathbf{f} \rangle}$ [with $\mathbf{f} = (f_1, \ldots, f_d)$] as the image measure—defined on $(\mathbb{R}, \mathscr{B}(\mathbb{R}))$—of $\alpha$, through the measurable function $\langle \mathbf{t}, \mathbf{f} \rangle := \sum_{j=1}^d t_j f_j$, that is,

$$\alpha_{\langle \mathbf{t}, \mathbf{f} \rangle} = \alpha \circ \langle \mathbf{t}, \mathbf{f} \rangle^{-1}.$$

If $\tilde{m}$ is a random probability measure on $(\Omega, \mathscr{F})$ with distribution $\mathscr{D}_\alpha$, then Theorem 9 states that the p.d. of

$$\tilde{m}(\langle \mathbf{t}, \mathbf{f} \rangle) = \int_\Omega \sum_{j=1}^d t_j f_j \, d\tilde{m}$$

coincides with the p.d. of $\tilde{\varphi}(I)$, where $\tilde{\varphi}$ is a random probability measure on $(\mathbb{R}, \mathscr{B}(\mathbb{R}))$ having p.d. $\mathscr{D}_{\alpha_{\langle \mathbf{t}, \mathbf{f} \rangle}}$.

8.1. *Multidimensional Lauricella identity.* These remarks, combined with a straightforward application of Theorem 4, yields the following *multidimensional form of the Lauricella identity*, that is: *Suppose* (8.1) *holds true, with* $a = \alpha(\Omega) > 0$, *and denote the p.d. of* $(\tilde{m}(f_1), \ldots, \tilde{m}(f_d))$ *by* $\mu_{\alpha, \mathbf{f}}$. *Then, for any* $\mathbf{t} = (t_1, \ldots, t_d)$ *in* $\mathbb{R}^d$, *one has*

$$\int_{\mathbb{R}^d} \frac{1}{(1 + i\langle \mathbf{t}, \mathbf{x} \rangle)^c} \mu_{\alpha, \mathbf{f}}(\, d\mathbf{x})$$

$$= \int_{[0,1]} \exp \left\{ - \int_\Omega \log(1 + iu\langle \mathbf{t}, \mathbf{f} \rangle) \, d\alpha \right\} B(\, du; c, a - c)$$

*if* $a \geq c > 0$, *and*

$$\int_{\mathbb{R}^d} \frac{1}{(1 + i\langle \mathbf{t}, \mathbf{x} \rangle)^c} \mu_{\alpha, \mathbf{f}}(\, d\mathbf{x})$$

$$= \frac{\Gamma(c - a + 1)\Gamma(a)}{2\pi i \Gamma(c)}$$

$$\times \int_0^{(1+)} \exp \left\{ - \int_\Omega \log(1 + iw\langle \mathbf{t}, \mathbf{f} \rangle) \, d\alpha \right\} w^{c-1}(w-1)^{a-c-1} \, dw$$

*if* $c > a > 0$.

This proposition has also been proved in Regazzini, Guglielmi and Di Nunno (2002) when $a = c$ and in Kerov and Tsilevich (1998) when $a = c = 1$. The following sections illustrate two applications. The former concerns the proof of absolute continuity, with respect to the Lebesgue measure on $\mathbb{R}^d$, of $\mu_{\alpha,\mathbf{f}}$. Regarding this point, Firmani (2002) has determined a density function of $\mu_{\alpha,\mathbf{f}}$ by inversion of the above identities with $d = 2$, and for suitable choices of $f_1$ and $f_2$. The latter application deals with the problem of determining an expression for the moment generating function of the variance of a Dirichlet process.

8.2. *Absolute continuity of $\mu_{\alpha,\mathbf{f}}$.* In view of Proposition 2 and of the remarks at the beginning of this section, the p.d. of $\tilde{m}(\langle \mathbf{t}, \mathbf{f} \rangle)$ is absolutely continuous with respect to the Lebesgue measure on $\mathbb{R}$ if $\mathbf{f} := (f_1, \ldots, f_d)$ is not *affinely $\alpha$-degenerate*; that is to say, there are no $\mathbf{v}$ in $\mathbb{R}^d \setminus \{0\}$ and $b$ in $\mathbb{R}$ for which $\alpha\{\langle \mathbf{v}, \mathbf{f} \rangle = b\} = a$. As far as absolute continuity of $\mu_{\alpha,\mathbf{f}}$ is concerned, the following statement is valid.

THEOREM 10. *Suppose $\alpha$ and $\mathbf{f}$ satisfy* (8.1) *and assume that $\mathbf{f}$ is not affinely $\alpha$-degenerate. Then the p.d. $\mu_{\alpha,\mathbf{f}}$ of $(\tilde{m}(f_1), \ldots, \tilde{m}(f_d))$ is absolutely continuous with respect to the Lebesgue measure on $\mathbb{R}^d$.*

PROOF. The contrapositive statement will be proved. Indeed, if $\mu_{\alpha,\mathbf{f}}$ is not absolutely continuous, there are $\varepsilon > 0$ and, for every $\delta > 0$, a set $A = \bigcup_{k=1}^{M} \Delta_k$ with $\Delta_i = \times_{r=1}^{d}(a_r^{(i)}, b_r^{(i)})$, $\Delta_i \cap \Delta_j = \varnothing$ if $i \neq j$, and $\lambda^d(A) = \sum_{k=1}^{M} \prod_{r=1}^{d}(b_r^{(k)} - a_r^{(k)}) < \delta$, such that $\mathscr{D}_\alpha\{\tilde{m}(\mathbf{f}) \in A\} = \mu_{\alpha,\mathbf{f}}(A) \geq \varepsilon$. Denote the Lebesgue measure on $\mathbb{R}^d$ by $\lambda^d$ and choose $t_1, \ldots, t_d$ in such a way that $0 < t_j \leq \min_{1 \leq k \leq M} \prod_{i \in \{1,\ldots,j-1,j+1,\ldots,d\}} (b_k^{(i)} - a_k^{(i)})/d$. Then,

$$\sum_{k=1}^{M} \sum_{j=1}^{d} t_j (b_k^{(j)} - a_k^{(j)}) \leq \sum_{k=1}^{M} \sum_{j=1}^{d} (b_k^{(j)} - a_k^{(j)}) < \delta$$

and, therefore, letting $A' = \bigcup_{k=1}^{M}(t_1 a_k^{(1)} + \cdots + t_d a_k^{(d)}, t_1 b_k^{(1)} + \cdots + t_d b_k^{(d)})$, one has $\lambda^1(A') \leq \sum_{k=1}^{M} \sum_{j=1}^{d} t_j (b_k^{(j)} - a_k^{(j)}) < \delta$. Moreover, $A \subset \{x \in \mathbb{R}^d : \langle \mathbf{t}, \mathbf{x} \rangle \in A'\}$ and $\varepsilon > 0$ is such that, with each $\delta > 0$, one can associate a set $A' \in \mathscr{B}(\mathbb{R})$ satisfying $\lambda^1(A') < \delta$ and $\mathscr{D}_\alpha\{\tilde{m}(\langle \mathbf{t}, \mathbf{f} \rangle) \in A'\} \geq \varepsilon$, contradicting the absolute continuity of the p.d. of $\tilde{m}(\langle \mathbf{t}, \mathbf{f} \rangle)$. $\quad\square$

8.3. *Moment generating function of the variance of $\tilde{\varphi}$.* If $\alpha$ satisfies condition (ii) in Theorem 1, then the random variance

$$\tilde{V} := \int_{\mathbb{R}} (x - \tilde{\varphi}(I))^2 \tilde{\varphi}(dx)$$

is finite, a.s.-$\mathscr{D}_\alpha$. Some of the results given in Section 7 are employed here to determine an expression for the moment generating function, $g_{\tilde{V}}$, of $\tilde{V}$.

The definition used for $g_{\tilde{V}}$ coincides with the one formulated, for instance, in Section 13.5 in Fristedt and Gray (1997), so that

$$(8.2) \quad g_{\tilde{V}}(t) = \mathscr{D}_\alpha(e^{-t\tilde{V}}) = \mathscr{D}_\alpha(\exp\{-t(\tilde{\varphi}(f_2) - (\tilde{\varphi}(f_1))^2)\}), \qquad t \geq 0,$$

where $f_1(x) = x$ and $f_2(x) = x^2$ for every $x$ in $\mathbb{R}$. After observing that $e^{tx^2}$ is the moment generating function, at $x$, of a Gaussian-distributed random variable with zero mean and variance equal to $2t$, apply Fubini's theorem to obtain

$$g_{\tilde{V}}(t) = \frac{1}{2\sqrt{2\pi t}} \int_{\mathbb{R}} e^{-u^2/(4t)} \mathscr{D}_\alpha(e^{u\varphi(f_1) - t\varphi(f_2)}) \, du.$$

If the support of $\alpha$ is bounded, by using arguments similar to those employed in the proof of Theorem 6 and by resorting to Theorem 9, with $f = uf_1 - tf_2$, it is straightforward to prove that

$$\mathscr{D}_\alpha(e^{u\varphi(f_1) - t\varphi(f_2)})$$
$$= \frac{\Gamma(a)}{2\pi i} \mathrm{PV} \int_{\gamma_{u,t} - i\infty}^{\gamma_{u,t} + i\infty} \frac{e^z}{z^a} \exp\left\{ -\int \log\left(1 - \frac{ux - tx^2}{z}\right) \alpha(dx) \right\} dz$$

holds, provided that $\gamma_{u,t}$ is any element of $(0, +\infty)$ satisfying

$$\gamma_{u,t} > \sup_{x \in \mathrm{supp}(\alpha)} (ux - tx^2).$$

The resulting expression for the moment generating function of $\tilde{V}$ is

$$\frac{\Gamma(a)}{4\pi^{3/2} i} \int_{\mathbb{R}} e^{-y^2/4} \mathrm{PV} \int_{\gamma'_{y,t} - i\infty}^{\gamma'_{y,t} + i\infty} \frac{e^z}{z^a} \exp\left\{ -\zeta\left(\frac{1}{z}; \alpha, tf_2 - y\sqrt{t} f_1\right) \right\} dz \, dy$$

and $\gamma'_{y,t}$ is any element of $(0, +\infty)$ satisfying

$$\gamma'_{y,t} > \sup_{x \in \mathrm{supp}(\alpha)} (y\sqrt{t} x - tx^2).$$

An extension of this representation to cases in which the support of $\alpha$ is arbitrary can be obtained by considering the sequence of truncated measures $(\alpha^{(k)})_{k \geq 1}$, as defined in Section 7. Since the distribution of the variance of the Dirichlet process with parameter $\alpha^{(k)}$ converges weakly to the distribution of the variance of the Dirichlet process with parameter $\alpha$, as $k$ tends to $+\infty$, a continuity theorem for moment generating functions yields

$$g_{\tilde{V}}(t) = \frac{\Gamma(a)}{4\pi^{3/2} i}$$
$$\times \lim_{k \to +\infty} \int_{\mathbb{R}} e^{-y^2/4}$$
$$\times \mathrm{PV} \int_{\gamma'_{y,t} - i\infty}^{\gamma'_{y,t} + i\infty} \frac{e^z}{z^a} \exp\left\{ -\zeta\left(\frac{1}{z}; \alpha^{(k)}, tf_2 - y\sqrt{t} f_1\right) \right\} dz \, dy.$$

See, for example, Fristedt and Gray [(1997), page 262].

**Acknowledgments.** Part of the present paper was written while the two authors were visiting Stanford University. They are grateful to people in the Department of Statistics for their kind hospitality and assistance. The authors also wish to thank two anonymous referees for their valuable comments and criticisms which have led to an improvement of the presentation.

DIPARTIMENTO DI ECONOMIA POLITICA
E METODI QUANTITATIVI
UNIVERSITÀ DEGLI STUDI DI PAVIA
VIA SAN FELICE 5
27100 PAVIA
ITALY
E-MAIL: lijoi@unipv.it

DIPARTIMENTO DI MATEMATICA
UNIVERSITÀ DEGLI STUDI DI PAVIA
VIA FERRATA 1
27100 PAVIA
ITALY
E-MAIL: eugenio@dimat.unipv.it